# MOMENTS OF MINORS OF WISHART MATRICES

By Mathias Drton,[1] Hélène Massam[2] and Ingram Olkin[3]

*University of Chicago, York University and Stanford University*

For a random matrix following a Wishart distribution, we derive formulas for the expectation and the covariance matrix of compound matrices. The compound matrix of order $m$ is populated by all $m \times m$-minors of the Wishart matrix. Our results yield first and second moments of the minors of the sample covariance matrix for multivariate normal observations. This work is motivated by the fact that such minors arise in the expression of constraints on the covariance matrix in many classical multivariate problems.

**1. Introduction.** Conditional independence constitutes one of the key concepts in multivariate statistical modeling. In a multivariate normal random vector $X = (X_1, \ldots, X_r)^T \sim \mathcal{N}_r(\mu, \Sigma)$, conditional independence expresses itself in the vanishing of minors, that is, subdeterminants of the positive definite covariance matrix. Let $I, J, K \subseteq [r] := \{1, \ldots, r\}$ be three pairwise disjoint index sets. Then $X_I$ and $X_J$ are conditionally independent given $X_K$, in symbols $X_I \perp\!\!\!\perp X_J \mid X_K$, if and only if

(1.1) $\qquad \det(\Sigma_{\{i\} \cup K \times \{j\} \cup K}) = 0 \qquad \forall i \in I, j \in J.$

The restrictions (1.1) correspond to vanishing partial correlations and can thus be tested using sample partial correlations, which yields a simple approach to model selection and assessment of goodness of fit of Gaussian independence models.

The situation becomes more complicated, however, in hidden variable models because conditional independences involving hidden variables may lead to constraints on the covariance matrix of the observed variables that no

Received June 2007; revised June 2007.
[1]Supported by NSF Grant DMS-05-05612.
[2]Supported by NSERC Discovery Grant A8946.
[3]Supported in part by NSF Grant REC-0634016.

longer correspond to vanishing partial correlations. Spearman's [14] *tetrads* are the classic example of such constraints. A tetrad is a $2 \times 2$-minor $\det(\Sigma_{ij \times k\ell})$ for which $\{i,j\} \cap \{k,\ell\} = \varnothing$. Tetrads are the defining equality constraints for one-factor analysis [6], but also arise in other Gaussian hidden variable models. (Recall that in factor analysis, observed variables are conditionally independent given hidden factors.)

Given a sample from a $\mathcal{N}_r(\mu, \Sigma)$ distribution (joint) vanishing of tetrads can be tested. Rejection of this hypothesis suggests that the model for which the tetrads would vanish is inappropriate for the data. The route commonly taken when testing the vanishing of a tetrad is to standardize the sample tetrad and compare the result to the standard normal distribution. This approach allows one in particular to avoid numerical maximization of the complicated likelihood functions of hidden variable models and we refer the reader to the examples discussed, for example, in [3, 8, 15]. The difficulty in this procedure is how to standardize the sample tetrad, a problem solved by Wishart [16] who found the sampling variance of the tetrad.

However, Wishart's result only applies to $2 \times 2$-minors, which has limited the application of the above constraint-based inference approach. In this paper we greatly generalize Wishart's result to obtain the covariance matrix of higher-order minors of a Wishart matrix, a problem that is also of intrinsic distribution-theoretic interest. In Section 2, we clarify the role of higher-order minors in hidden variable models. In Section 3, we present some basic results based on simple but powerful invariance arguments for compound matrices. Together with the properties of the Choleski decomposition of a Wishart matrix, these results allow us to compute, in Sections 4 and 5, the expectations and covariance matrix of minors of arbitrary Wishart matrices. In our conclusion in Section 6, we comment on future research directions and give an example of constraint-based inference based on $3 \times 3$-minors. The Appendix contains the proofs of two lemmas as well as an interesting auxiliary result on the mean of the determinant of a noncentral Wishart covariance matrix.

**2. Off-diagonal minors and hidden variables.** Tetrads are $2 \times 2$-minors that do not involve any diagonal elements of the covariance matrix $\Sigma$. We call any minor with this property an *off-diagonal* minor. In seminal work, Spirtes, Glymour and Scheines [15], Theorem 6.10, have characterized the tetrad relations in covariance matrices from directed Gaussian graphical models. The characterization of the vanishing of higher-order off-diagonal minors is still an open problem but in Proposition 2.2 below we are able to give simple sufficient conditions. Proposition 2.2(ii) applies in particular to factor analysis with $m - 1$ factors; see also [6].

Consider a random vector $X = (X_1, \ldots, X_r)^T \sim \mathcal{N}_r(\mu, \Sigma)$ with $r \geq 2m$ components. Let $I, J \subseteq [r]$ be two disjoint index sets of cardinality $|I| = |J| = m \geq 1$.

LEMMA 2.1. *If $K \subseteq [r] \setminus (I \cup J)$, then $X_I \perp\!\!\!\perp X_J \mid X_K$ if and only if $\det(\Sigma_{G \times H}) = 0$ for all $G \subseteq I \cup K$ and $H \subseteq J \cup K$ of cardinality $|G| = |H| = |K| + 1$.*

PROOF. The claimed vanishing of minors implies (1.1) and thus the conditional independence. Conversely, the conditional independence implies that

$$\Sigma_{I \cup K \times K \cup J} = \begin{pmatrix} \Sigma_{I \times K} & \Sigma_{I \times J} \\ \Sigma_{K \times K} & \Sigma_{K \times J} \end{pmatrix} = \begin{pmatrix} \Sigma_{I \times K} \\ \Sigma_{K \times K} \end{pmatrix} \Sigma_{K \times K}^{-1} \begin{pmatrix} \Sigma_{K \times K} & \Sigma_{K \times J} \end{pmatrix}.$$

The second equality implies $\mathrm{rank}(\Sigma_{I \cup K \times J \cup K}) \leq |K|$ and thus the claim. □

PROPOSITION 2.2. (i) *If $X_i \perp\!\!\!\perp X_J$ for some $i \in I$, then $\det(\Sigma_{I \times J}) = 0$.*
(ii) *Suppose the partitions $I = I_1 \dot\cup I_2$ and $J = J_1 \dot\cup J_2$ have $I_1, J_1 \neq \varnothing$ ($I_2$ or $J_2$ may be empty). Let $K_1$ and $K$ be subsets of $[r] \setminus (I \cup J)$ such that $K_1 \subseteq K$ and $|K| + |I_2| + |J_2| \leq |K_1| + m - 1$. If $X_{I_1} \perp\!\!\!\perp X_{J_1} \mid X_{K \cup I_2 \cup J_2}$ and $X_I \perp\!\!\!\perp X_{K_1}$, then $\det(\Sigma_{I \times J}) = 0$. The proposition states in particular that if $X_I \perp\!\!\!\perp X_J \mid X_K$ for $K \subseteq [r] \setminus (I \cup J)$ with $|K| \leq m-1$, then $\det(\Sigma_{I \times J}) = 0$.*

PROOF. (i) Immediate. (ii) By Lemma 2.1, $\mathrm{rank}(\Sigma_{I \cup K \times J \cup K}) \leq |K| + |I_2| + |J_2|$ and thus $\det(\Sigma_{I \cup K_1 \times J \cup K_1}) = 0$. Since $\Sigma_{I \times K_1} = 0$, it holds that $\det(\Sigma_{I \cup K_1 \times J \cup K_1}) = \det(\Sigma_{I \times J}) \det(\Sigma_{K_1 \times K_1})$, and the claim follows because $\det(\Sigma_{K_1 \times K_1}) > 0$. The last statement of the proposition is obtained from (ii) by taking $I_2 = J_2 = K_1 = \varnothing$. □

For an example in which an $m \times m$-minor yields the only equality constraint on the covariance matrix, consider $2m + (m-1)$ random variables $X_1, \ldots, X_{2m}, Y_1, \ldots, Y_{m-1}$. Define an acyclic digraph (DAG) $G_m$ with these random variables as vertices and edges as follows. Every variable $Y_i$ is adjacent to every one of the variables $X_j$ by a directed edge $Y_i \to X_j$. Every pair of vertices in $\{X_1, \ldots, X_m\}$ is joined by an edge, and the same holds for every pair of vertices in $\{X_{m+1}, \ldots, X_{2m}\}$. For uniqueness assume that $X_i \to X_j$ implies $i < j$. Figure 1 shows the graph $G_3$, which we will take up in a data example in the conclusion in Section 6.

In the remainder of this section, let $I = [m]$ and $J = \{m+1, \ldots, 2m\}$. The graph $G_m$ encodes that $X_I$ is conditionally independent of $X_J$ given $Y_{[m-1]}$ and that the random variables $Y_i$ are completely independent; see, for example, [10] for details on the stochastic interpretation of directed graphs. Treating $Y_1, \ldots, Y_{m-1}$ as hidden yields a Gaussian model for $(X_1, \ldots, X_{2m})^T$. It can be shown that this model contains exactly those distributions $\mathcal{N}_{2m}(\mu, \Sigma)$ that have a covariance matrix of the form $\Sigma = \Omega + \Lambda\Lambda^T$, where $\Lambda$ is an arbitrary $2m \times (m-1)$-matrix and $\Omega$ is a positive definite block-diagonal matrix; $\Omega_{I \times J} = 0$. Let $\mathcal{C}_m$ be the set of covariance matrices $\Sigma$ in this model. The following lemma is proven in the Appendix.

LEMMA 2.3. *If a polynomial $f$ in the entries of the covariance matrix $\Sigma$ evaluates to zero at every matrix in $\mathcal{C}_m$, then $f$ is a polynomial multiple of the off-diagonal $m \times m$-minor $\det(\Sigma_{I \times J})$.*

**3. Invariance under orthogonal transformations.** Given that minors of covariance matrices arise so naturally in independence models, it is interesting to study their natural estimators, namely the minors of the sample covariance matrix. Up to a scaling factor depending on the sample size, such sample minors are distributed like the minors of Wishart matrices, which arise as follows.

Let $X \in \mathbb{R}^{r \times n}$ be a matrix whose columns are independent random vectors distributed according to the multivariate normal distribution $\mathcal{N}_r(0, \Sigma)$ with positive definite covariance matrix $\Sigma \in \mathbb{R}^{r \times r}$. Then $S = XX^T$ is distributed according to the *Wishart distribution* with scale parameter matrix $\Sigma$ and $n$ degrees of freedom, in symbols, $S \sim \mathcal{W}_r(n, \Sigma)$. We refer to the Wishart distribution $\mathcal{W}_r(n, I_r)$ with the identity matrix $I_r \in \mathbb{R}^{r \times r}$ as scale parameter, as *standard Wishart distribution*.

Simple invariance arguments based on ideas from Olkin and Rubin [13] (see also [4] and [7], Problem 4, page 330) will permit us to learn much about the standard case.

DEFINITION 3.1. *Let $O(r)$ be the group of orthogonal matrices in $\mathbb{R}^{r \times r}$. The distribution of a symmetric random matrix $V \in \mathbb{R}^{r \times r}$ is orthogonally invariant, if for all $G \in O(r)$, the distribution of $GVG^T$ is identical to the distribution of $V$. We will say, for brevity, that $V \in \mathbb{R}^{r \times r}$ is orthogonally invariant.*

For $S \sim \mathcal{W}_r(m, \Sigma)$ and $G \in O(r)$, we have that $GSG^T \sim \mathcal{W}_r(n, G\Sigma G^T)$, and hence, the standard Wishart distribution $\mathcal{W}_r(n, I_r)$ is orthogonally invariant.

The objects of our study are minors $\det(W_{I \times J})$ or $\det(S_{I \times J})$ that are specified by two subsets $I, J \subseteq [r]$ of equal cardinality $|I| = |J| = m$. We

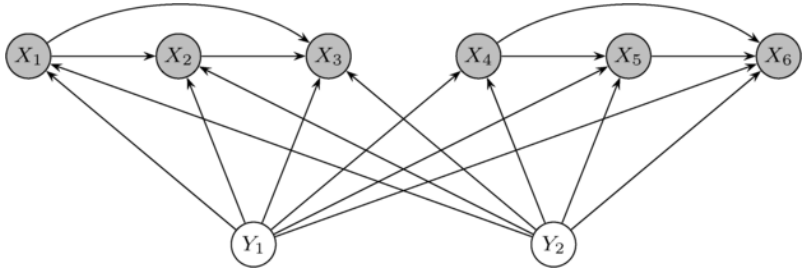

FIG. 1. *The graph $G_3$ with two complete subgraphs joined through two hidden variables.*

introduce the notation

$$\left\{{r \atop m}\right\} = \{I \subseteq [r] : |I| = m\}, \qquad m \in [r].$$

PROPOSITION 3.2. *Let $I, J \in \left\{{r \atop m}\right\}$. If the symmetric random matrix $V \in \mathbb{R}^{r \times r}$ is orthogonally invariant, then*

$$\mathrm{E}[\det(V_{I \times J})] = \begin{cases} \mathrm{E}[\det(V_{[m] \times [m]})], & \text{if } I = J, \\ 0, & \text{otherwise.} \end{cases}$$

PROOF. We extend the proof of [13], Lemma 1, which treats the case $m = 1$, in which the minors reduce to individual entries of $V$. Let $I, J \in \left\{{r \atop m}\right\}$ be two distinct subsets. For $j \in J \setminus I$, let $D_j \in O(r)$ be the diagonal matrix equal to the identity matrix except for entry $(j, j)$ which is equal to $-1$. Then $D_j V D_j^T$ differs from $V$ in that all off-diagonal entries of the $j$th row and column have been negated. Since $j \in J$ but $j \notin I$, $\det[(D_j V D_j^T)_{I \times J}] = -\det(V_{I \times J})$. Thus, $\mathrm{E}[\det(V_{I \times J})] = 0$ is implied by $\mathrm{E}[\det(V_{I \times J})] = \mathrm{E}[\det(D_j V D_j^T)_{I \times J}] = -\mathrm{E}[\det(V_{I \times J})]$.

Since $|I| = |J|$, we can find a permutation that maps the indices in $I$ to those in $J$. Let $P = P_{IJ} \in O(r)$ be the matrix representing this permutation. Then, $(PVP^T)_{I \times I} = V_{J \times J}$, and $\mathrm{E}[\det(V_{I \times I})] = \mathrm{E}\{\det[(PVP^T)_{I \times I}]\} = \mathrm{E}[\det(V_{J \times J})]$. It follows that $\mathrm{E}[\det(V_{I \times I})] = \mathrm{E}[\det(V_{[m] \times [m]})]$ for all $I \in \left\{{r \atop m}\right\}$. □

Our approach to determine the moment structure of the minors of a Wishart matrix is based on the following ideas. First, recall that if $S \sim \mathcal{W}_r(n, \Sigma)$, and $W \sim \mathcal{W}_r(n, I_r)$, then

(3.1)  $$S = \Sigma^{1/2} W \Sigma^{1/2} \sim \mathcal{W}_r(n, \Sigma^{1/2} I_r \Sigma^{1/2}) = \mathcal{W}_r(n, \Sigma).$$

[For notational simplicity, we will use the symmetric square root throughout the paper but nonsymmetric square roots (e.g., lower triangular) could be used instead.] Second, recall that for a matrix $A \in \mathbb{R}^{r \times r}$ and an integer $m \in [r]$, the $m$th compound of $A$ is the matrix

$$A^{(m)} = (\det(A_{I \times J}))_{I, J \in \left\{{r \atop m}\right\}} \in \mathbb{R}^{\binom{p}{m} \times \binom{p}{m}}$$

that is populated with all $m \times m$-minors of $A$. If $m = 0$, we set $A^{(0)} = 1 \in \mathbb{R}$. The Binet–Cauchy theorem (see, e.g., Marshall and Olkin [11], page 503 and Aitken [1], Chapter V) states that

(3.2)  $$(AB)^{(m)} = A^{(m)} B^{(m)},$$

which allows us to use (3.1) for the transfer from standard to general Wishart matrices. The last ingredient to our approach is the fact that the products

$$\det(S_{I\times J})\det(S_{K\times L}), \qquad I,J,K,L \in \left\{{r \atop m}\right\},$$

which are exactly the quantities of interest for studying the variance-covariance structure of minors of $S$, are the entries of the Kronecker product $S^{(m)} \otimes S^{(m)}$.

The next proposition states that the first and second moments of the compound matrix $S^{(m)}$ can be obtained from those of the compound matrix $W^{(m)}$ for a standard Wishart matrix. The result follows from (3.1) and (3.2).

PROPOSITION 3.3. *Let $S \sim \mathcal{W}_r(n,\Sigma)$ and $W \sim \mathcal{W}_r(n,I_r)$ and let $\Sigma^{1/2}$ denote the unique symmetric and positive definite square root of $\Sigma$. Then*

$$\mathrm{E}[S^{(m)}] = (\Sigma^{1/2})^{(m)}\mathrm{E}[W^{(m)}](\Sigma^{1/2})^{(m)},$$

$$\mathrm{E}[S^{(m)} \otimes S^{(m)}] = [(\Sigma^{1/2})^{(m)} \otimes (\Sigma^{1/2})^{(m)}]$$
$$\times (\mathrm{E}[W^{(m)} \otimes W^{(m)}])[(\Sigma^{1/2})^{(m)} \otimes (\Sigma^{1/2})^{(m)}],$$

$$\mathrm{Cov}[S^{(m)}] = \mathrm{E}[S^{(m)} \otimes S^{(m)}] - (\mathrm{E}[S^{(m)}] \otimes \mathrm{E}[S^{(m)}])$$
$$= [(\Sigma^{1/2})^{(m)} \otimes (\Sigma^{1/2})^{(m)}](\mathrm{Cov}[W^{(m)}])[(\Sigma^{1/2})^{(m)} \otimes (\Sigma^{1/2})^{(m)}].$$

Proposition 3.3 is only useful if we are able to compute the necessary moments of $W^{(m)}$. However, the invariance of $W$ under the orthogonal group tells us a great deal about these moments. The full first and second moment structure of $W^{(m)}$ will be derived in Corollary 4.2 and Theorem 4.5.

In the next result $I \triangle J$ denotes the symmetric difference $(J \setminus I) \cup (I \setminus J)$.

PROPOSITION 3.4. *Let $I, J, K, L \in \{{r \atop m}\}$, and let $V \in \mathbb{R}^{r \times r}$ be orthogonally invariant. If $I \triangle J \neq K \triangle L$, then*

(3.3) $$\mathrm{E}[\det(V_{I\times J})\det(V_{K\times L})] = 0.$$

*Moreover, under any permutation of $\sigma$ the indices in $[r]$,*

(3.4) $$\mathrm{E}[\det(V_{I\times J})\det(V_{K\times L})] = \mathrm{E}[\det(V_{\sigma(I)\times\sigma(J)})\det(V_{\sigma(K)\times\sigma(L)})].$$

PROOF. Again we extend the ideas in [13], Lemma 1.

Let $I \triangle J \neq K \triangle L$. Assume without loss of generality that there exists an index $j \in (I \triangle J) \setminus (K \triangle L)$. Let $D_j$ be the diagonal matrix defined in the proof of Proposition 3.2. Recall that the action of $D_j$ negates the $j$th row and column in $V$. By choice of $j \in I \triangle J$, it holds that $\det[(D_j V D_j^T)_{I \times J}] = -\det(V_{I \times J})$. Since either $j \in K \cap L$ or $j \notin K \cup L$, it holds further that $\det[(D_j V D_j^T)_{K \times L}] = \det(V_{K \times L})$. Then (3.3), as well as (3.4), follows immediately from the orthogonal invariance of $V$. $\square$

EXAMPLE 3.5. Let $m = 2$ and $r = 4$. Applying the permutation $\sigma = (1)(23)(4)$, (3.4) implies that $\mathrm{E}[\det(V_{12 \times 12})^2] = \mathrm{E}[\det(V_{13 \times 13})^2]$ but because there is no permutation of $\{1,2,3,4\}$ that sends $\{1,2\}$ to $\{1,2\}$ and $\{1,2\}$ to $\{2,4\}$, $\mathrm{E}[\det(V_{12 \times 12})^2]$ and $\mathrm{E}[\det(V_{12 \times 12}) \det(V_{24 \times 24})]$ need not be equal. We will illustrate the use of Proposition 3.4 further in Example 4.6 where we consider a standard Wishart distribution.

**4. Choleski-decomposition of a standard Wishart matrix.** The arguments presented in Section 3 determine the first and second moments of minors of orthogonally invariant random matrices only up to constants. In this section we determine these constants for the standard Wishart distribution $\mathcal{W}_r(n, I_r)$. For this task we use the Choleski-decomposition that has the following convenient distributional property; see, for example, Muirhead [12], Theorem 3.2.14, for a proof of this classical result.

LEMMA 4.1. *Let $W$ follow the $\mathcal{W}_r(n, I_r)$ distribution with $n \geq r$. Let $T = (t_{ij})_{i \leq i, j \leq r}$ be lower-triangular with positive diagonal entries such that $W = TT^T$. Then the $t_{ij}$, $i \geq j$, are independent random variables distributed as $t_{ii}^2 \sim \chi_{n-i+1}^2, i = 1, \ldots, r$, and $t_{ij} \sim \mathcal{N}(0,1), 1 \leq j < i \leq r$.*

We remark that the elements $t_{ij}$ have been called rectangular coordinates. Since $\det(TT^T) = \prod_{i=1}^r t_{ii}^2$ and $\mathrm{E}[t_{ii}^2] = \mathrm{E}[\chi_{n-i+1}^2] = n - i + 1$, we obtain the following corollary to Proposition 3.2 and to Lemma 4.1 applied to $W_{I \times I}$.

COROLLARY 4.2. *Let $I, J \in \{^r_m\}$. If $W \sim \mathcal{W}_r(n, I_r)$ with $n \geq m$, then*

$$\mathrm{E}[\det(W_{I \times J})] = \begin{cases} n!/(n-m)!, & \text{if } I = J, \\ 0, & \text{otherwise,} \end{cases}$$

*and*

$$\mathrm{E}[S^{(m)}] = \frac{n!}{(n-m)!} \Sigma^{(m)}.$$

The Choleski-decomposition $W = TT^T$ of a standard Wishart matrix $W$ reveals additional information. In the remainder of this section assume that $n \geq r$, which implies that $T$ is of full rank with probability 1.

LEMMA 4.3. *Let $c \in \{0, \ldots, m\}$ be an integer such that there exists a subset $\bar{J} \subseteq \{m+1, \ldots, r\}$ of cardinality $|\bar{J}| = m - c$. Then*

$$\det(W_{[m] \times ([c] \cup \bar{J})}) = \left(\prod_{i=1}^{c} t_{ii}^2\right)\left(\prod_{j=c+1}^{m} t_{jj}\right) \det(T_{\bar{J} \times \{c+1,\ldots,m\}}).$$

PROOF. Let $\bar{I} = \{c+1, \ldots, m\} = [m] \setminus [c]$. From the partitioning

$$W_{[m] \times ([c] \cup \bar{J})} = \begin{pmatrix} W_{[c] \times [c]} & W_{[c] \times \bar{J}} \\ W_{\bar{I} \times [c]} & W_{\bar{I} \times \bar{J}} \end{pmatrix},$$

we obtain that

(4.1) $\det(W_{[m] \times ([c] \cup \bar{J})}) = \det(W_{[c] \times [c]}) \det(W_{\bar{I} \times \bar{J}} - W_{\bar{I} \times [c]} W_{[c] \times [c]}^{-1} W_{[c] \times \bar{J}}).$

Clearly, $\det(W_{[c] \times [c]}) = \prod_{i=1}^{c} t_{ii}^2$ so that we are left with studying the second factor on the right-hand side of (4.1).

For a subset $D \subseteq [r]$, let $T_D = T_{D \times [r]}$ be the submatrix comprising all rows of $T$ with index in $D$. Then we can write

$$\det(W_{\bar{I} \times \bar{J}} - W_{\bar{I} \times [c]} W_{[c] \times [c]}^{-1} W_{[c] \times \bar{J}}) = \det\{T_{\bar{I}}[I_r - T_{[c]}^T (T_{[c]} T_{[c]}^T)^{-1} T_{[c]}] T_{\bar{J}}^T\}.$$

The matrix $I_r - T_{[c]}^T (T_{[c]} T_{[c]}^T)^{-1} T_{[c]}$ represents the orthogonal projection on the kernel of $T_{[c]}$. Since $T$ is lower diagonal with by assumption nonzero diagonal entries, it holds that $\ker(T_{[c]}) = \{0\}^c \times \mathbb{R}^{r-c}$, which means that the projection considered replaces the first $c$ entries of a vector in $\mathbb{R}^r$ by zeros. Therefore,

$$\det\{T_{\bar{I}}[I_r - T_{[c]}^T (T_{[c]} T_{[c]}^T)^{-1} T_{[c]}] T_{\bar{J}}^T\} = \det(T_{\bar{I} \times \{c+1,\ldots,r\}} T_{\bar{J} \times \{c+1,\ldots,r\}}^T).$$

By the Binet–Cauchy theorem,

$$\det(T_{\bar{I} \times \{c+1,\ldots,r\}} T_{\bar{J} \times \{c+1,\ldots,r\}}^T) = \sum_{\substack{D \subseteq \{c+1,\ldots,r\}, \\ |D|=m-c}} \det(T_{\bar{I} \times D}) \det(T_{\bar{J} \times D})$$

(4.2)
$$= \det(T_{\bar{I} \times \bar{I}}) \det(T_{\bar{J} \times \bar{I}}).$$

The second equality in (4.2) holds because if $D \neq \bar{I} = \{c+1, \ldots, m\}$, then the matrix $T_{\bar{I} \times D}$ contains a column consisting entirely of zeros and thus $\det(T_{\bar{I} \times D}) = 0$. Our claim follows from $\det(T_{\bar{I} \times \bar{I}}) = \prod_{j=c+1}^{m} t_{jj}$. □

From Lemma 4.3 we can deduce the distribution of a minor $\det(W_{I \times J})$.

THEOREM 4.4. *Let $I, J \in \binom{r}{m}$ have a (possibly empty) intersection of cardinality $|I \cap J| = c \geq 0$. If $W \sim \mathcal{W}_r(n, I_r)$, then*

(4.3) $$\det(W_{I \times J}) \sim \left(\prod_{i=1}^{c} W_i\right)\left(\prod_{i=c+1}^{m} \sqrt{W_i}\right) \det(Z),$$

where $W_i, i = 1, \ldots, m$, are independent $\chi^2_{n-i+1}$ random variables and $Z = (Z_{ij}) \in \mathbb{R}^{(m-c) \times (m-c)}$ is a random matrix of independent $\mathcal{N}(0,1)$ random variables that are also independent of $(W_1, \ldots, W_m)$. In particular,

$$\mathrm{E}[\det(W_{I \times J})^2] = \frac{n!}{(n-m)!} \frac{(n+2)!}{(n+2-|I \cap J|)!} (m - |I \cap J|)!$$

and

$$\mathrm{Var}[\det(W_{I \times J})]$$
$$= \begin{cases} \dfrac{n!}{(n-m)!} \left[ \dfrac{(n+2)!}{(n+2-m)!} - \dfrac{n!}{(n-m)!} \right], & \text{if } I = J, \\ \dfrac{n!}{(n-m)!} \dfrac{(n+2)!}{(n+2-|I \cap J|)!} (m - |I \cap J|)!, & \text{if } |I \cap J| < m. \end{cases}$$

PROOF. By orthogonal invariance of the standard Wishart distribution, we can permute rows and columns of $W$ such that $I = [m]$ and $J = [c] \cup \{m+1, \ldots, 2m-c\}$. Thus (4.3) follows from Lemmas 4.3 and 4.1.

For the derivation of the second moment, recall that $\mathrm{E}[\chi^2_n] = n$ and $\mathrm{E}[(\chi^2_n)^2] = n(n+2)$. Let $\mathcal{S}_h$ be the group of permutations of $[h]$. Then

$$\mathrm{E}[\det(Z)^2] = \sum_{\sigma \in \mathcal{S}_{m-c}} \sum_{\tau \in \mathcal{S}_{m-c}} \mathrm{E}\left[ \prod_{i=1}^{m-c} Z_{i\sigma(i)} Z_{i\tau(i)} \right]$$
$$= \sum_{\sigma \in \mathcal{S}_{m-c}} \prod_{i=1}^{m-c} \mathrm{E}[Z^2_{i\sigma(i)}] = (m-c)!,$$

which yields $\mathrm{E}[\det(W_{I \times J})^2]$. The variance is obtained using Corollary 4.2. □

We next turn to moments of the form $\mathrm{E}[\det(W_{I \times J}) \det(W_{K \times L})]$ with $(I, J) \neq (K, L)$. By Proposition 3.4, this expectation is nonzero only if $I \triangle J = K \triangle L$. At this point, before proceeding to derive the desired expectation, we would like to emphasize that throughout the paper we consider an index set $I = \{i_1, \ldots, i_m\}$ to be equipped with an ordering. Such an ordering yields an index sequence $(i_1, \ldots, i_m)$ that dictates the order in which we list the rows (or columns) of a submatrix. Since our results so far did not depend on the choice of ordering, we kept this view implicit. For our next result, however, the order in which the indices in $I$ are listed matters since different orderings may lead to different signs of determinants due to the interchanging of the rows or columns in submatrices. For example,

(4.4) $\quad \mathrm{E}[\det(W_{12 \times 14}) \det(W_{23 \times 34})] = -\mathrm{E}[\det(W_{12 \times 14}) \det(W_{23 \times 43})].$

In the following theorem, the elements of four index sets $I$, $J$, $K$, $L$ are assumed to be ordered according to a total order of $[r]$ that achieves certain order relationships across the four sets. We write $A < B$ if all elements of $A \subseteq [r]$ are smaller than those of $B \subseteq [r]$, or if $A$ or $B$ is the empty set. (Note that this implies $\varnothing < A < \varnothing$.)

THEOREM 4.5. *Let $I, J, K, L \in \{{}^r_m\}$ such that $I \triangle J = K \triangle L$. Let*
$$\bar{I} = I \setminus (I \cap J), \qquad \bar{K} = K \setminus (K \cap L), \qquad \bar{J} = J \setminus (I \cap J), \qquad \bar{L} = L \setminus (K \cap L).$$
*Moreover, assume that the indices in $I$, $J$, $K$ and $L$ are listed according to a total order in $[r]$ under which*
$$(I \cap J) \setminus (K \cap L) < \bar{I} < \bar{J} < (K \cap L) \setminus (I \cap J),$$
$$\bar{I} \cap \bar{K} < \bar{I} \cap \bar{L}, \qquad \bar{J} \cap \bar{K} < \bar{J} \cap \bar{L}.$$
*Under these conventions it holds that if $W \sim \mathcal{W}_r(n, I_r)$, then*
$$\mathrm{E}[\det(W_{I \times J}) \det(W_{K \times L})] = \frac{n!}{(n-m)!} \frac{(n+2)!}{(n+2-|I \cap J \cap K \cap L|)!}$$
$$\times \frac{(n-m+|(I \cap J) \setminus (K \cap L)|)!}{(n-m)!} |\bar{I} \cap \bar{K}|! |\bar{I} \cap \bar{L}|!.$$

Theorem 4.5 yields, for example, that $\mathrm{E}[\det(W_{12 \times 13}) \det(W_{24 \times 34})] = n(n-1)^2$. However, it does not yield directly the value of $\mathrm{E}[\det(W_{12 \times 14}) \times \det(W_{23 \times 34})]$ in (4.4). Instead, we can obtain that $\mathrm{E}[\det(W_{12 \times 14}) \det(W_{23 \times 43})] = n(n-1)^2$. Hence, by (4.4), we find $\mathrm{E}[\det(W_{12 \times 14}) \det(W_{23 \times 34})] = -n(n-1)^2$.

EXAMPLE 4.6. If $m = 2$ and $r = 4$, then the covariance matrix $\mathrm{Cov}[W^{(m)}]$, which determines $\mathrm{Cov}[S^{(m)}]$, is a symmetric matrix of size $36 \times 36$. Because $\mathrm{Cov}[W^{(m)}]$ is derived from the symmetric matrix $W$, we can restrict ourselves to unordered pairs of sets $(I, J) \in \{{}^r_m\} \times \{{}^r_m\}$ with possible equality $I = J$. There are 21 such unordered pairs. We represent $\mathrm{Cov}[W^{(m)}]$ as a symmetric block-diagonal $21 \times 21$ matrix with blocks formed according to $I \triangle J$. Proposition 3.4 implies such block-diagonal structure also for the general case of arbitrary $m$ and $r$.

The first block is indexed by the six pairs $(I, I)$, $I \in \{{}^r_m\}$, involves the principal minors and takes on the form

$$\begin{pmatrix} 12,12 & 13,13 & 14,14 & 23,23 & 24,24 & 34,34 \\ f_1 & f_4 & f_4 & f_4 & f_4 & 0 \\ & f_1 & f_4 & f_4 & 0 & f_4 \\ & & f_1 & 0 & f_4 & f_4 \\ & & & f_1 & f_4 & f_4 \\ & & & & f_1 & f_4 \\ & & & & & f_1 \end{pmatrix},$$

where from Theorems 4.4 and 4.5 respectively, we have $f_1 = 2n(2n+1)(n-1)$ and $f_4 = 2n(n-1)^2$. Next, we have a series of six blocks of size $2 \times 2$, each involving two pairs $(I, J)$ and $(K, L)$ for which $I \triangle J = K \triangle L$ and $|I \cap J| = 1$, or equivalently, $|I \triangle J| = 2$. Two representatives of these six blocks are

$$\begin{array}{cc} 12,13 & 24,34 \\ \begin{pmatrix} f_2 & f_5 \\ & f_2 \end{pmatrix} \end{array} \quad \text{and} \quad \begin{array}{cc} 12,14 & 23,34 \\ \begin{pmatrix} f_2 & -f_5 \\ & f_2 \end{pmatrix}, \end{array}$$

where by Theorems 4.4 and 4.5 respectively, $f_2 = n(n+2)(n-1)$ and $f_5 = n(n-1)^2$. The last block is obtained for the pairs $(I, J)$ with $I, J$ disjoint, or equivalently, $I \triangle J = [r] = \{1, 2, 3, 4\}$. It takes the form

$$\begin{array}{ccc} 12,34 & 13,24 & 14,23 \\ \begin{pmatrix} f_3 & f_6 & -f_6 \\ & f_3 & f_6 \\ & & f_3 \end{pmatrix} \end{array}$$

with $f_3 = 2n(n-1)$ and $f_6 = n(n-1)$.

The remainder of this section is devoted to the proof of Theorem 4.5 in which we can assume that $r = \max(I \cup J \cup K \cup L)$. Note that since $|(I \cap J) \setminus (K \cap L)| = |(K \cap L) \setminus (I \cap J)|$, the formula in Theorem 4.5 is not changed if the order of $(I, J)$ and $(K, L)$ is reversed.

LEMMA 4.7. *If* $I \cap J \cap K \cap L = C \neq \varnothing$, $|C| = c \geq 1$, *then*

$$\mathrm{E}[\det(W_{I \times J}) \det(W_{K \times L})] = \mathrm{E}[\det(\bar{W}_{I^c \times J^c}) \det(\bar{W}_{K^c \times L^c})] \cdot \mathrm{E}[\det(W_{C \times C})^2],$$

*where* $A^c = A \setminus C$ *for any subset* $A \subseteq [r]$, *and*

$$\bar{W} = W_{[r]^c \times [r]^c} - W_{[r]^c \times C} W_{C \times C}^{-1} W_{C \times [r]^c} \sim \mathcal{W}_{r-c}(n - c, I_{r-c}).$$

PROOF. The claim follows from the fact that

$$\det(W_{I \times J}) \det(W_{K \times L}) = \det(W_{C \times C})^2 \det(\bar{W}_{I^c \times J^c}) \det(\bar{W}_{K^c \times L^c})$$

in conjunction with the independence of $W_{C \times C}$ and $\bar{W}$ (see Lemma 5.2 below). □

Since Theorem 4.4 yields the term $\mathrm{E}[\det(W_{C \times C})^2]$ appearing in Lemma 4.7, the proof of Theorem 4.5 is completed by the following lemma, which is proven in the Appendix.

LEMMA 4.8. *Let* $I, J, K, L \in \binom{r}{m}$ *such that* $I \triangle J = K \triangle L$ *and* $I \cap J \cap K \cap L = \varnothing$. *Define* $\bar{I}, \bar{J}, \bar{K}, \bar{L}$ *as in Theorem 4.5, and assume furthermore* that $I \cap J < \bar{I} < \bar{J} < K \cap L$, $\bar{I} \cap \bar{K} < \bar{I} \cap \bar{L}$, and $\bar{J} \cap \bar{K} < \bar{J} \cap \bar{L}$. If $W \sim \mathcal{W}_r(n, I_r)$, then

$$\mathrm{E}[\det(W_{I \times J}) \det(W_{K \times L})] = \frac{n!(n-m+c)!}{[(n-m)!]^2} \cdot p!(m-c-p)!,$$

where $c = |I \cap J| = |K \cap L|$ and $p = |\bar{I} \cap \bar{K}| = |\bar{J} \cap \bar{L}|$.

**5. Variances of minors.** In Sections 3–4 we found the covariance matrix of the compound $S^{(m)}$ of a Wishart matrix $S \sim \mathcal{W}_r(n, \Sigma)$. However, due to the involved square roots $\Sigma^{1/2}$, the form of the individual entries of $\mathrm{Cov}[S^{(m)}]$ is not transparent. In this section, we derive explicit formulas for the variances of $m \times m$-minors with $m \leq n$.

We begin by reviewing the well-known formula for a principal minor [2], Section 7.5.

PROPOSITION 5.1. *If $S \sim \mathcal{W}_r(n, \Sigma)$ and $I \in \binom{r}{m}$, then*

$$\mathrm{Var}[\det(S_{I \times I})] = \frac{n!}{(n-m)!} \left\{ \frac{(n+2)!}{(n+2-m)!} - \frac{n!}{(n-m)!} \right\} \det(\Sigma_{I \times I})^2.$$

PROOF. Apply (3.1) with the submatrix $S_{I \times I}$ replacing the full Wishart matrix $S_{I \times I}$ to obtain that $\mathrm{Var}[\det(S_{I \times I})] = \det(\Sigma_{I \times I})^2 \cdot \mathrm{Var}[\det(W_{I \times I})]$, which in conjunction with Theorem 4.4 yields the claim. □

Next, we derive an explicit formula for the variance of off-diagonal minors of a general Wishart matrix $S \sim \mathcal{W}_r(n, \Sigma)$. From this formula and Proposition 5.1, a formula for the variance of arbitrary minors of $S$ is obtained in Theorem 5.7.

Let $I, J \in \binom{r}{m}$ be two *disjoint* subsets. Then the minor $\det(S_{I \times J})$ is off-diagonal in that it does not involve any diagonal elements of $S$. Let $S_{IJ \times IJ}$ and $\Sigma_{IJ \times IJ}$ be the $(I \cup J) \times (I \cup J)$-submatrix of $S$ and $\Sigma$, respectively. We partition these $2m \times 2m$-submatrices into four $m \times m$-submatrices according to $I$ and $J$ where we adopt the shorthand notation $S_{I \times I} = S_{II}$, $S_{I \times J} = S_{IJ}$, etc. Let

$$S_{II.J} = S_{II} - S_{IJ} S_{JJ}^{-1} S_{JI} \quad \text{and} \quad \Sigma_{II.J} = \Sigma_{II} - \Sigma_{IJ} \Sigma_{JJ}^{-1} \Sigma_{JI}.$$

Our line of attack in computing the variance of the off-diagonal minor $\det(S_{I \times J}) = \det(S_{IJ})$ is to employ the decomposition

(5.1) $\mathrm{Var}[\det(S_{IJ})] = \mathrm{Var}[\mathrm{E}[\det(S_{IJ}) \mid S_{JJ}]] + \mathrm{E}[\mathrm{Var}[\det(S_{IJ}) \mid S_{JJ}]].$

The evaluations of the two terms on the right-hand side of (5.1) are given in Lemmas 5.3 and 5.4, which are based on the following well-known result [12], Theorem 3.2.10.

LEMMA 5.2. *If $S \sim \mathcal{W}_r(n, \Sigma)$ and $m < n$ then $S_{JJ} \sim \mathcal{W}_m(n, \Sigma_{JJ})$, $S_{II.J} \sim \mathcal{W}_m(n-m, \Sigma_{II.J})$, and the random matrix $S_{II.J}$ is independent of $(S_{IJ}, S_{JJ})$. Finally, the conditional distribution of $S_{IJ}$ given $S_{JJ}$ is normal and such that*

$$(5.2) \qquad (S_{IJ} S_{JJ}^{-1/2} \mid S_{JJ}) \sim \mathcal{N}_{m^2}(\Sigma_{IJ} \Sigma_{JJ}^{-1} S_{JJ}^{1/2}, \Sigma_{II.J} \otimes I_m)$$

$$(5.3) \iff (\Sigma_{II.J}^{-1/2} S_{IJ} S_{JJ}^{-1/2} \mid S_{JJ}) \sim \mathcal{N}_{m^2}(\Sigma_{II.J}^{-1/2} \Sigma_{IJ} \Sigma_{JJ}^{-1} S_{JJ}^{1/2}, I_m \otimes I_m).$$

LEMMA 5.3. *It holds that*

$$\operatorname{Var}[\operatorname{E}[\det(S_{IJ}) \mid S_{JJ}]] = \frac{n!}{(n-m)!} \left\{ \frac{(n+2)!}{(n+2-m)!} - \frac{n!}{(n-m)!} \right\} \cdot \det(\Sigma_{IJ})^2.$$

PROOF. By (5.2) in Lemma 5.2,

$$\operatorname{Var}[\operatorname{E}[\det(S_{IJ}) \mid S_{JJ}]] = \operatorname{Var}[\operatorname{E}[\det(S_{IJ} S_{JJ}^{-1/2}) \mid S_{JJ}] \cdot \det(S_{JJ}^{1/2})]$$
$$= \operatorname{Var}[\det(\Sigma_{IJ} \Sigma_{JJ}^{-1}) \cdot \det(S_{JJ})]$$
$$= \det(\Sigma_{IJ})^2 \det(\Sigma_{JJ})^{-2} \operatorname{Var}[\det(S_{JJ})].$$

Now the claim follows from Proposition 5.1. □

LEMMA 5.4. *Let $\Sigma^{IJ}$ denote the $I \times J$-submatrix of the inverse of $\Sigma_{IJ \times IJ}$. Then*

$$\operatorname{E}[\operatorname{Var}[\det(S_{IJ}) \mid S_{JJ}]]$$
$$= \det(\Sigma_{IJ \times IJ})$$
$$\times \left( \sum_{k=0}^{m-1} (m-k)! \cdot \frac{n!}{(n-m)!} \cdot \frac{(n+2)!}{(n+2-k)!} \cdot (-1)^k \operatorname{tr}\{(\Sigma_{JI} \Sigma^{IJ})^{(k)}\} \right).$$

PROOF. First note that

$$(5.4) \quad \begin{aligned} &\operatorname{E}[\operatorname{Var}[\det(S_{IJ}) \mid S_{JJ}]] \\ &= \det(\Sigma_{II.J}) \cdot \operatorname{E}[\operatorname{Var}[\det(\Sigma_{II.J}^{-1/2} S_{IJ} S_{JJ}^{-1/2}) \mid S_{JJ}] \cdot \det(S_{JJ})]. \end{aligned}$$

It follows from (5.3) that conditional on $S_{JJ}$, the entries of the matrix $\Sigma_{II.J}^{-1/2} S_{IJ} S_{JJ}^{-1/2}$ are independent normal random variables with variance 1, albeit these entries are not identically distributed as their means may differ in arbitrary fashion. We are led to the problem of computing $\operatorname{Var}[\det(X)]$, where the matrix $X \in \mathbb{R}^{m \times m}$ is distributed according to the multivariate normal distribution

$$X \sim \mathcal{N}_{m^2}(A, I_m \otimes I_m), \qquad A = \Sigma_{II.J}^{-1/2} \Sigma_{IJ} \Sigma_{JJ}^{-1} S_{JJ}^{1/2} \in \mathbb{R}^{m \times m}.$$

Lemma A.1 provides an evaluation of $\mathrm{E}[\det(X)^2]$, and from (5.4) we find that

$$
\begin{aligned}
\mathrm{E}[\mathrm{Var}[\det(S_{IJ}) \mid S_{JJ}]] \\
= \det(\Sigma_{II.J}) \\
\times \sum_{k=0}^{m-1} (m-k)! \\
\times \mathrm{E}\left[\mathrm{tr}\{(\Sigma_{II.J}^{-1/2}\Sigma_{IJ}\Sigma_{JJ}^{-1}S_{JJ}\Sigma_{JJ}^{-1}\Sigma_{JI}\Sigma_{II.J}^{-1/2})^{(k)}\} \cdot \det(S_{JJ})\right].
\end{aligned}
\tag{5.5}
$$

Setting $C = \Sigma_{JJ}^{-1}\Sigma_{JI}\Sigma_{II.J}^{-1}\Sigma_{IJ}\Sigma_{JJ}^{-1}$, (5.5) simplifies to

$$
\begin{aligned}
\mathrm{E}[\mathrm{Var}[\det(S_{IJ}) \mid S_{JJ}]] &= \det(\Sigma_{II.J}) \sum_{k=0}^{m-1} (m-k)! \cdot \mathrm{E}[\mathrm{tr}\{(CS_{JJ})^{(k)}\} \cdot \det(S_{JJ})] \\
&= \det(\Sigma_{II.J}) \sum_{k=0}^{m-1} (m-k)! \cdot \mathrm{E}[\mathrm{tr}\{C^{(k)} S_{JJ}^{(k)} \cdot \det(S_{JJ})\}] \\
&= \det(\Sigma_{II.J}) \sum_{k=0}^{m-1} (m-k)! \cdot \mathrm{tr}\{C^{(k)} \mathrm{E}[S_{JJ}^{(k)} \cdot \det(S_{JJ})]\}.
\end{aligned}
$$

Now, let $W_{JJ} = (\Sigma_{JJ})^{-1/2} S_{JJ} (\Sigma_{JJ})^{-1/2}$. As in (3.1), $W_{JJ} \sim \mathcal{W}_m(n, I_m)$. Thus

$$
\begin{aligned}
&\mathrm{E}[\mathrm{Var}[\det(S_{IJ}) \mid S_{JJ}]] \\
&= \det(\Sigma_{II.J}) \det(\Sigma_{JJ}) \\
&\quad \times \left( \sum_{k=0}^{m-1} (m-k)! \cdot \mathrm{tr}\{C^{(k)} (\Sigma_{JJ}^{1/2})^{(k)} \mathrm{E}[W_{JJ}^{(k)} \cdot \det(W_{JJ})] (\Sigma_{JJ}^{1/2})^{(k)}\} \right).
\end{aligned}
$$

The distribution of $W_{JJ}^{(k)} \cdot \det(W_{JJ})$ has the invariance property that for $G \in O(m)$,

$$G^{(k)}(W_{JJ}^{(k)} \cdot \det(W_{JJ}))(G^T)^{(k)} \sim W_{JJ}^{(k)} \cdot \det(W_{JJ}).$$

In analogy to Proposition 3.2 and the derivation of Theorem 4.4, it holds that

$$\mathrm{E}[W_{JJ}^{(k)} \cdot \det(W_{JJ})] = \frac{n!}{(n-m)!} \cdot \frac{(n+2)!}{(n+2-k)!} \cdot I_{\binom{m}{k}}.$$

Because $\det(\Sigma_{II.J})\det(\Sigma_{JJ}) = \det(\Sigma_{IJ \times IJ})$, we therefore have that

$\mathrm{E}[\mathrm{Var}[\det(S_{IJ}) \mid S_{JJ}]]$

$$= \det(\Sigma_{IJ \times IJ}) \cdot \left( \sum_{k=0}^{m-1} (m-k)! \frac{n!}{(n-m)!} \cdot \frac{(n+2)!}{(n+2-k)!} \cdot \operatorname{tr}\{(\Sigma_{JJ}C)^{(k)}\} \right).$$

The claim now follows because, by simple considerations about the inverse of the partitioned matrix $\Sigma_{IJ \times IJ}$, it holds that $\Sigma_{JJ}C = -\Sigma_{JI}\Sigma^{IJ}$. $\square$

Combining Lemmas 5.3 and 5.4 according to (5.1) yields the following proposition.

PROPOSITION 5.5. *Let $I, J \in \{^r_m\}$ be two disjoint subsets. Then the off-diagonal minor $\det(S_{I \times J}) = \det(S_{IJ})$ of the Wishart matrix $S \sim \mathcal{W}_r(n, \Sigma)$ has variance*

$$\operatorname{Var}[\det(S_{IJ})] = \frac{n!}{(n-m)!} \cdot \det(\Sigma_{IJ})^2 \left\{ \frac{(n+2)!}{(n+2-m)!} - \frac{n!}{(n-m)!} \right\}$$

$$+ \frac{n!}{(n-m)!} \cdot \det(\Sigma_{IJ \times IJ})$$

$$\times \left( \sum_{k=0}^{m-1} (m-k)! \cdot \frac{(n+2)!}{(n+2-k)!} \cdot (-1)^k \operatorname{tr}\{(\Sigma_{JI}\Sigma^{IJ})^{(k)}\} \right).$$

COROLLARY 5.6 ([16]). *In the special case $m=2$ the off-diagonal minor $\det(S_{I \times J}) = \det(S_{IJ})$ is known as a tetrad, and*

$$\operatorname{Var}[\det(S_{IJ})] = n(n-1)[(n+2)\det(\Sigma_{II})\det(\Sigma_{JJ})$$
$$- n \det(\Sigma_{IJ \times IJ}) + 3n \det(\Sigma_{IJ})^2].$$

PROOF. The claim follows from Proposition 5.5, and the fact that if $m = 2$, then

$$\operatorname{tr}(\Sigma_{JI}\Sigma^{IJ}) \det(\Sigma_{IJ \times IJ}) = \det(\Sigma_{IJ \times IJ}) - \det(\Sigma_{II})\det(\Sigma_{JJ}) + \det(\Sigma_{IJ})^2. \quad \square$$

THEOREM 5.7. *Let $I, J \in \{^r_m\}$ have intersection $C := I \cap J$ of cardinality $c = |C| = |I \cap J|$. Define $\bar{I} = I \setminus (I \cap J)$, $\bar{J} = J \setminus (I \cap J)$ and $\bar{I}\bar{J} = \bar{I} \cup \bar{J}$. Then the minor $\det(S_{I \times J}) = \det(S_{IJ})$ of the Wishart matrix $S \sim \mathcal{W}_r(n, \Sigma)$ has variance*

$$\operatorname{Var}[\det(S_{IJ})] = \frac{n!}{(n-m)!} \left\{ \frac{(n+2)!}{(n+2-c)!} - \frac{n!}{(n-c)!} \right\} \det(\Sigma_{C \times C})^2$$

$$\times \left[ \det(\bar{\Sigma}_{\bar{I}\bar{J}})^2 \left\{ \frac{(n+2-c)!}{(n+2-m)!} - \frac{(n-c)!}{(n-m)!} \right\} \right.$$

$$+ \det(\bar{\Sigma}_{\bar{I}\bar{J} \times \bar{I}\bar{J}})$$

$$\times \left( \sum_{k=0}^{m-c-1} (m-c-k)! \cdot \frac{(n+2-c)!}{(n+2-c-k)!} \cdot (-1)^k \right.$$
$$\left. \times \operatorname{tr}\{(\bar{\Sigma}_{\bar{J}\bar{I}} \bar{\Sigma}^{\bar{I}\bar{J}})^{(k)}\} \right) \Bigg],$$

where $\bar{\Sigma} = \Sigma_{([r]\backslash C) \times ([r]\backslash C)} - \Sigma_{([r]\backslash C) \times C} \Sigma_{C \times C}^{-1} \Sigma_{C \times ([r]\backslash C)}$.

PROOF. Define $\bar{S}$ in analogy to $\bar{\Sigma}$. Since $\det(S_{I \times J}) = \det(S_{C \times C}) \det(\bar{S}_{\bar{I} \times \bar{J}})$ and $S_{C \times C}$ and $\bar{S}_{\bar{I} \times \bar{J}}$ are independent (Lemma 5.2), the claim follows from Propositions 5.1 and 5.5. □

**6. Conclusion.** We study first and second moments of minors of a Wishart matrix, relying fundamentally on the properties of compound matrices. For a standard Wishart matrix $W$, invariant under $O(r)$, we extended classic invariance arguments due to Olkin and Rubin [13] to the case of compounds. This was possible because the Binet–Cauchy theorem implies that the distribution of the compound $W^{(m)}$ is invariant under compounds of matrices in $O(r)$. Note, however, that the distribution of $W^{(m)}$ is not invariant under all matrices in $O(\binom{r}{m})$.

Our results yield closed-form test statistics that are useful for evaluating the goodness of fit of hidden variable models; compare [6], Section 3. As an example, consider the model from Section 2 that is induced by the graph $G_3$ in Figure 1. For illustration we use classic data on physical variables for 305 fifteen-year-old girls from the University of Chicago Lab schools; a correlation matrix is reported in [9], Table 7.1, page 169. We choose $X_1, \ldots, X_6$ as *Height*, *Arm span*, *Length of forearm*, *Weight*, *Chest girth* and *Chest width*. The partition in $I = \{X_1, X_2, X_3\}$ and $J = \{X_4, X_5, X_6\}$ thus separates variables relating to lankiness from those relating to stockiness. We compute the $I \times J$ minor of the sample correlation matrix (recall Lemma 2.3) and estimate its sampling variance by inserting the correlation matrix into the formula from Proposition 5.5. When doing this we omit the first term in the formula because $\det(\Sigma_{IJ})$ is hypothesized to be zero. Comparing the ratio of sample minor and estimated standard deviation to the standard normal distribution gives a $p$-value of 0.42. In comparison, the likelihood ratio test computed using the EM algorithm has $p$-value 0.39, which also indicates a good model fit. Repeating the same procedure for a less meaningful variable partition obtained by exchanging *Length of forearm* ($X_3$) and *Chest width* ($X_6$) leads to $p$-values of 0.0034 and 0.0026 for the minor and the likelihood ratio test, respectively. These results suggest that the closed-form minor test may indeed have good power.

In the above example, the only data available were a sample correlation matrix, which we treated as if it were a sample covariance matrix. This is justified, however, because the ratio of sample minor and standard deviation estimate is the same when evaluated over the sample correlation matrix instead of the sample covariance matrix. This fact is a consequence of the multilinearity of the determinant and the Binet–Cauchy theorem, which implies that $\mathrm{Var}_\Sigma[\det(S_{I\times J})] = (\prod_{i\in I}\sigma_{ii})(\prod_{j\in J}\sigma_{jj})\mathrm{Var}_R[\det(S_{I\times J})]$. Here, $R$ is the correlation matrix of the covariance matrix $\Sigma$. While we can justifiably compute standardized sample minors from correlation matrices, our Wishart distribution results do not yield the moments of minors of sample correlation matrices. The determination of these is an interesting problem for future research. The distribution of sample correlation matrices is orthogonally invariant when the covariance matrix is a multiple of the identity but it is not so in general.

Our data example falls into a traditional large sample setting. We believe that minors may also be useful for high-dimensional settings in which the number of variables is large, perhaps even larger than the sample size. The reasoning behind this speculation is that sample minors may be formed from full rank submatrices even when the entire sample covariance matrix is singular. Clearly a likelihood ratio test against a saturated alternative is impossible under such singularity.

## APPENDIX

**A.1. Proof of Lemma 2.3.** Let $\mathbb{R}[\sigma]$ be the ring of polynomials in the indeterminates $\sigma_{ij}$, $i \leq j$. Define $\mathcal{I}_1 \subseteq \mathbb{R}[\sigma]$ to be the ideal generated by the minor $\det(\Sigma_{I\times J})$. Since this minor is irreducible, $\mathcal{I}_1$ is a prime ideal. Define $\mathcal{I}_2 \subseteq \mathbb{R}[\sigma]$ to be the ideal of all polynomials that vanish when evaluated at a matrix $\Sigma \in \mathcal{C}_m$. The ideal $\mathcal{I}_2$ is also a prime ideal [5], Section 4.5. In Lemma 2.3, we claim that $\mathcal{I}_1 = \mathcal{I}_2$.

Let $V_1$ and $V_2$ be the irreducible varieties of complex matrices $\Sigma$ such that $f(\Sigma) = 0$ for all $f \in \mathcal{I}_1$ and all $f \in \mathcal{I}_2$, respectively. In all distributions in the graphical model induced by the graph $G_m$ defined in Section 2 it holds that $X_I \perp\!\!\!\perp X_J \mid Y_{[m-1]}$. Hence, by Proposition 2.2(ii), $\det(\Sigma_{I\times J}) = 0$ for all $\Sigma \in \mathcal{C}_m$, which implies that $\mathcal{I}_1 \subseteq \mathcal{I}_2$ and $V_2 \subseteq V_1$. Conversely, a matrix $\Sigma \in V_1$ can be written as $\Sigma = \Omega + \Lambda\Lambda^T$ with $\Omega \in \mathbb{C}^{2m\times 2m}$ block-diagonal and $\Lambda \in \mathbb{C}^{2m\times(m-1)}$. A polynomial in $\mathcal{I}_2$ must vanish at such a matrix $\Sigma$. Thus $\Sigma \in V_2$, and consequently $V_1 = V_2$. Since $\mathcal{I}_1$ is a prime ideal it now follows from the Strong Nullstellensatz [5], Section 4.2, that $\mathcal{I}_1 = \mathcal{I}_2$.

**A.2. Proof of Lemma 4.8.** First, we emphasize that $\bar{I} \cap \bar{J} = \varnothing$, $\bar{K} \cap \bar{L} = \varnothing$, $\bar{I} \dot\cup \bar{J} = I \triangle J = K \triangle L = \bar{K} \dot\cup \bar{L}$, and $|\bar{I}| = |\bar{J}| = |\bar{K}| = |\bar{L}| = m-c$. Defining $q = |\bar{I} \cap \bar{L}| = |\bar{J} \cap \bar{K}|$, it also holds that $p + q = |\bar{I} \cap \bar{K}| + |\bar{I} \cap \bar{L}| = |\bar{I}| = m-c$. Moreover, since $|\bar{I} \cap \bar{K}| + |\bar{J} \cap \bar{K}| = |\bar{K}| = m-c$, it holds that $|\bar{I} \cap \bar{K}| = p = |\bar{J} \cap \bar{L}|$ and $|\bar{I} \cap \bar{L}| = q = m-c-p = |\bar{J} \cap \bar{K}|$.

By permuting the indices in $[r]$ if necessary (Proposition 3.4), we can assume that

$$I \cap J = \{1, \ldots, c\},$$
$$\bar{I} \cap \bar{K} = \{c+1, \ldots, c+p\},$$
$$\bar{I} \cap \bar{L} = \{c+p+1, \ldots, m = c+p+q\},$$
$$\bar{J} \cap \bar{K} = \{m+1, \ldots, m+q+1\},$$
$$\bar{J} \cap \bar{L} = \{m+q+1, \ldots, 2m-c = m+q+p\},$$
$$K \cap L = \{2m-c+1, \ldots, 2m\}.$$

As another convention, we enumerate the elements of the sets $K$ and $L$ as $K = (k_1, \ldots, k_m)$ and $L = (\ell_1, \ldots, \ell_m)$, respectively, while choosing $k_i = \ell_i$ for all $i \in [c]$.

Let $W = TT^T$ be the Choleski-decomposition of $W$ whose Choleski-factor $T = (t_{ij})$ is lower-triangular with positive diagonal elements. By Lemma 4.3,

$$(A.1) \qquad \det(W_{I \times J}) = \left(\prod_{i=1}^{c} t_{ii}^2\right)\left(\prod_{i=c+1}^{m} t_{ii}\right) \det(T_{\bar{J} \times \bar{I}}).$$

Whereas $\det(W_{I \times J})$ has the simple representation in (A.1), this is not the case for $\det(W_{K \times L})$. However, because we are interested in $\mathrm{E}[\det(W_{I \times J}) \det(W_{K \times L})]$ some simplification is possible based on the following fact. Because $t_{ij}$, $i > j$, are independent $\mathcal{N}(0,1)$, if $(\alpha_{ij} \mid 1 \le j \le i \le r)$ contains an entry $\alpha_{ij}$ that is odd and such that $i > j$, then

$$(A.2) \qquad \mathrm{E}\left[\prod_{i \ge j} t_{ij}^{\alpha_{ij}}\right] = \prod_{i \ge j} \mathrm{E}[t_{ij}^{\alpha_{ij}}] = 0.$$

By the Binet–Cauchy theorem, $\det(W_{K \times L})$ is equal to

$$\sum_{H \in \{^r_m\}} \det(T_{K \times H}) \det(T_{L \times H}) = \sum_{H \in \{^r_m\}} \sum_{\sigma \in \mathcal{S}_m} \sum_{\tau \in \mathcal{S}_m} (-1)^{\sigma+\tau} \prod_{a=1}^{m} t_{k_a h_{\sigma(a)}} t_{\ell_a h_{\tau(a)}},$$

where $H = \{h_1, \ldots, h_m\}$. Since $k_a = \ell_a$ for $a \in [c]$,

$$(A.3) \qquad \prod_{a=1}^{c} t_{k_a h_{\sigma(a)}} t_{\ell_a h_{\tau(a)}} = \prod_{a=1}^{c} t_{k_a h_{\sigma(a)}} t_{k_a h_{\tau(a)}}.$$

We claim that

$$(A.4) \qquad f_1 = \sum_{H \in \{^r_m\}} \sum_{\sigma \in \mathcal{S}_m} \sum_{\substack{\tau \in \mathcal{S}_m: \\ \tau(a) = \sigma(a) \forall a \in [c]}} (-1)^{\sigma+\tau}$$
$$\times \left(\prod_{a=1}^{c} t_{k_a h_{\sigma(a)}}^2\right)\left(\prod_{b=c+1}^{m} t_{k_b h_{\sigma(b)}} t_{\ell_b h_{\tau(b)}}\right)$$

satisfies

$$\mathrm{E}[\det(W_{I\times J})\det(W_{K\times L})] = \mathrm{E}[\det(W_{I\times J})\cdot f_1].$$

In order to see this, fix $H$ and $\sigma$, and assume that $\tau$ is such that there exists $a \in [c]$ for which $\sigma(a) \neq \tau(a)$. Then $h_{\sigma(a)} \neq k_a$ or $h_{\tau(a)} \neq k_a$. Without loss of generality, assume that $h_{\sigma(a)} \neq k_a$. If $h_{\sigma(a)} > k_a$, then $t_{k_a h_{\sigma(a)}} = 0$ because $T$ is lower-triangular. If $h_{\sigma(a)} < k_a$, then $t_{k_a h_{\sigma(a)}}$ appears with exponent 1 in the monomial $\prod_{a=1}^m t_{k_a h_{\sigma(a)}} t_{\ell_a h_{\tau(a)}}$. The index $k_a \in K \cap L$ is not an element of $I \cup J$. Thus $t_{k_a h_{\sigma(a)}}$ appears with exponent 1 in

$$\det(W_{I\times J}) \cdot \prod_{a=1}^m t_{k_a h_{\sigma(a)}} t_{\ell_a h_{\tau(a)}}.$$

Therefore, according to (A.2), only monomials $\prod_{a=1}^m t_{k_a h_{\sigma(a)}} t_{\ell_a h_{\tau(a)}}$ appearing in $f_1$ may contribute to the expected value of $\det(W_{I\times J})\det(W_{K\times L})$.

We can rewrite (A.4) as

$$f_1 = \sum_{H \in \{{r \atop m}\}} \sum_{\sigma \in \mathcal{S}_m} (-1)^\sigma \left( \prod_{a=1}^c t_{k_a h_{\sigma(a)}}^2 \right) \left( \prod_{b=c+1}^m t_{k_b h_{\sigma(b)}} \right)$$

$$\times \left[ \sum_{\substack{\tau \in \mathcal{S}_m: \\ \tau(a)=\sigma(a) \forall a \in [c]}} (-1)^\tau \prod_{b=c+1}^m t_{\ell_b h_{\tau(b)}} \right].$$

Now,

$$\sum_{\substack{\tau \in \mathcal{S}_m: \\ \tau(a)=\sigma(a) \forall a \in [c]}} (-1)^\tau \left( \prod_{b=c+1}^m t_{\ell_b h_{\tau(b)}} \right)$$

$$= (-1)^\sigma \sum_{\substack{\tau \in \mathcal{S}_m: \\ \tau(a)=\sigma(a) \forall a \in [c]}} (-1)^{\tau \circ \sigma^{-1}} \prod_{b=c+1}^m t_{\ell_b h_{\tau \circ \sigma^{-1}(\sigma(b))}}$$

$$= (-1)^\sigma \sum_{\bar\tau \in \mathcal{S}_{\sigma(\{c+1,\ldots,m\})}} (-1)^{\bar\tau} \prod_{b=c+1}^m t_{\ell_b h_{\bar\tau(\sigma(b))}}$$

$$= (-1)^\sigma \det(T_{\bar L \times (h_{\sigma(c+1)},\ldots,h_{\sigma(m)})}).$$

Therefore,

$$f_1 = \sum_{H \in \{{r \atop m}\}} \sum_{\sigma \in \mathcal{S}_m} \left( \prod_{a=1}^c t_{k_a h_{\sigma(a)}}^2 \right) \left( \prod_{b=c+1}^m t_{k_b h_{\sigma(b)}} \right) \det(T_{L \times (h_{\sigma(c+1)},\ldots,h_{\sigma(m)})}).$$

We have thus shown that $\mathrm{E}[\det(W_{I\times J})\det(W_{K\times L})]$ is equal to the expectation of

$$\sum_{H\in\{^r_m\}}\sum_{\sigma\in\mathcal{S}_m}\left(\prod_{a=1}^{c}t^2_{k_a h_{\sigma(a)}}\right)\left(\prod_{b=c+1}^{m}t_{k_b h_{\sigma(b)}}\right)$$

(A.5)

$$\times\det(T_{L\times(h_{\sigma(c+1)},\ldots,h_{\sigma(m)})})\left(\prod_{a=1}^{c}t^2_{i_a i_a}\right)\left(\prod_{b=c+1}^{m}t_{i_b i_b}\right)\det(T_{\bar{J}\times\bar{I}}).$$

We next claim that the expectation of (A.5) does not change when dropping all terms associated with pairs $(H,\sigma)$ for which $\{h_{\sigma(c+1)},\ldots,h_{\sigma(m)}\}\neq\bar{I}$. To see this, choose $b\in\{c+1,\ldots,m\}$ for which $h_{\sigma(b)}\in\{h_{\sigma(c+1)},\ldots,h_{\sigma(m)}\}\setminus\bar{I}$. Now consider three cases. First, if $h_{\sigma(b)}\in(K\cap L)\,\dot{\cup}\,(\bar{J}\cap\bar{L})$, then $h_{\sigma(b)}>k_b\in\bar{K}$, and it follows that $t_{k_b h_{\sigma(b)}}=0$, which leads to the vanishing of the term associated with $H$ and $\sigma$. Second, if $h_{\sigma(b)}\in\bar{J}\cap\bar{K}$, then every nonzero term in the expansion of $\det(T_{L\times(h_{\sigma(c+1)},\ldots,h_{\sigma(m)})})$ involves an off-diagonal element of $T$ that does not appear in $\det(T_{\bar{J}\times\bar{I}})$. Hence, every monomial of the term associated with $(H,\sigma)$ features an off-diagonal element of $T$ raised to the power 1. Therefore, by (A.2), the term associated with $(H,\sigma)$ has expectation zero. The third case in which $h_{\sigma(b)}\in I\cap J$ is similar to the second case just discussed.

The claim just verified allows us to rewrite (A.5) as

$$\sum_{\substack{H\in\{^r_m\}:\bar{I}\subseteq H}}\sum_{\substack{\sigma\in\mathcal{S}_m,\\ h_{\sigma(\{c+1,\ldots,m\})}=\bar{I}}}(-1)^{\nu_\sigma}\left(\prod_{a=1}^{c}t^2_{k_a h_{\sigma(a)}}\right)$$

(A.6)
$$\times\left(\prod_{a=1}^{c}t^2_{i_a i_a}\right)\left(\prod_{b=c+1}^{m}t_{k_b h_{\sigma(b)}}\right)$$

$$\times\left(\prod_{b=c+1}^{m}t_{i_b i_b}\right)\det(T_{\bar{L}\times\bar{I}})\det(T_{\bar{J}\times\bar{I}}),$$

where $\nu_\sigma$ is the permutation of $\bar{I}=\{c+1,\ldots,m\}$ that sorts $h_{\sigma(c+1)},\ldots,h_{\sigma(m)}$ in increasing order, that is, $c+1=h_{\sigma(\nu(c+1))}<\cdots<h_{\sigma(\nu(m))}=m$. We now argue that the expectation of (A.6) does not change when replacing $\det(T_{\bar{L}\times\bar{I}})$ by

(A.7) $$\left(\prod_{\ell\in\bar{L}\cap\bar{I}}t_{\ell\ell}\right)\det(T_{(\bar{J}\cap\bar{L})\times(\bar{I}\cap\bar{K})}).$$

In fact, we find (A.7) from the Laplace expansion along the diagonal $t_{\ell\ell}$, $\ell\in\bar{L}\cap\bar{I}$, for which $\det(T_{(\bar{J}\cap\bar{L})\times(\bar{I}\times\bar{K})})$ serves as a cofactor. Every term in $\det(T_{\bar{L}\times\bar{I}})$ that does not appear in (A.7) involves an off-diagonal entry in $T$ of the form $t_{ab}$ with $a \in \bar{L} \cap \bar{I}$ and $b \in \bar{I}$, $a > b$. Such $t_{ab}$, however, does not appear in $\det(T_{\bar{J}\times\bar{I}})$ since clearly $a < \min(\bar{J})$. Now, an appeal to (A.2) closes the argument.

Next, recall that $k_b \in \bar{I} \cap \bar{K}$ if $b \in \{c+1, \ldots, c+p\}$. Hence, if $b \in \{c+1, \ldots, c+p\}$ but $h_{\sigma(b)} \neq k_b$, then the term $t_{k_b, h_{\sigma(b)}}$ does not appear in $\det(T_{\bar{J}\times\bar{I}})$. In other words, a term in (A.6) based on $(H, \sigma)$ with $h_{\sigma(b)} \neq k_b$ has zero expectation by (A.2). Combining this observation with the replacement in (A.7), we define

$$f_2 = \sum_{H \in \{^r_m\}: \bar{I} \subseteq H} \sum_{\substack{\sigma \in \mathcal{S}_m: h_{\sigma(\{c+1,\ldots,m\})} = \bar{I}, \\ h_{\sigma(c+j)} = k_{c+j} \in \bar{I} \cap \bar{K} \forall j \in [p]}} (-1)^{\nu_\sigma} \left(\prod_{a=1}^c t^2_{k_a h_{\sigma(a)}}\right) \left(\prod_{b=c+p+1}^m t_{k_b h_{\sigma(b)}}\right)$$

$$\text{(A.8)} \times \left(\prod_{k \in \bar{I} \cap \bar{K}} t_{kk}\right) \left(\prod_{a=1}^c t^2_{i_a i_a}\right) \left(\prod_{b=c+1}^m t_{i_b i_b}\right) \left(\prod_{\ell \in \bar{L} \cap \bar{I}} t_{\ell\ell}\right)$$

$$\times \det(T_{(\bar{J}\cap\bar{L})\times(\bar{I}\cap\bar{K})}) \det(T_{\bar{J}\times\bar{I}}),$$

which satisfies $E[f_2] = E[\det(W_{I\times J}) \det(W_{K\times L})]$.

In our next simplification, we claim that if we replace $\det(T_{\bar{J}\times\bar{I}})$ in $f_2$ by

$$\text{(A.9)} \qquad \det(T_{(\bar{J}\cap\bar{L})\times(\bar{I}\cap\bar{K})}) \det(T_{(\bar{J}\cap\bar{K})\times(\bar{I}\cap\bar{L})}),$$

then the expectation does not change. This follows from (A.2) because every term in the expansion of $\det(T_{\bar{J}\times\bar{I}})$ that does not appear in (A.9) involves some $t_{ab}$ with $a \in \bar{J} \cap \bar{L}$ and $b \in \bar{I} \cap \bar{L}$, and such $t_{ab}$ appears neither in $\det(T_{(\bar{J}\cap\bar{L})\times(\bar{I}\cap\bar{K})})$ nor in $\prod_{b=c+p+1}^m t_{k_b h_{\sigma(b)}}$ because $k_b \in \bar{K} < \min(\bar{J} \cap \bar{L})$.

In $f_2$, $H \in \{^r_m\}$ is such that $\bar{I} \subseteq H$ and $h_{\sigma(\{c+1,\ldots,m\})} = \bar{I}$ and therefore

$$\text{(A.10)} \qquad H = \{h_1, h_2, \ldots, h_c\} \cup \bar{I}.$$

Using (A.9) and the fact that $h_{\sigma(b)} \in \bar{I} \setminus (\bar{I} \cap \bar{K}) = \bar{I} \cap \bar{L}$ if $b \in \{c+p+1, \ldots, m\} = \bar{I} \cap \bar{L}$, we obtain that $E[\det(W_{I\times J}) \det(W_{K\times L})]$ is equal to

$$\left[\sum_{h_1 \in [k_1]\setminus\bar{I}} \sum_{h_2 \in [k_2]\setminus(\bar{I}\cup\{h_1\})} \cdots \sum_{h_c \in [k_c]\setminus(\bar{I}\cup\{h_1,\ldots,h_{c-1}\})} \left(\prod_{a=1}^c t^2_{k_a h_a}\right)\right]$$

$$\text{(A.11)} \qquad \times \left(\sum_{\mu \in \mathcal{S}_{\bar{I}\cap\bar{L}}} (-1)^\mu \prod_{b=c+p+1}^m t_{k_b \mu(b)}\right) \left(\prod_{a=1}^c t^2_{i_a i_a}\right) \left(\prod_{i \in \bar{I}} t^2_{ii}\right)$$

$$\times \det(T_{(\bar{J}\cap\bar{L})\times(\bar{I}\cap\bar{K})})^2 \det(T_{(\bar{J}\cap\bar{K})\times(\bar{I}\cap\bar{L})}).$$

In the simplification from (A.8) to (A.11) we replaced the two sums over $H$ and $\sigma$ by the sums over $h_1, \ldots, h_c$. This is possible because of (A.10) and because by choosing appropriate $H$ and $\sigma$, $h_{\sigma(a)}$ can take on any value in $[k_a] \setminus \bar{I}$ while respecting that all $h_{\sigma(a)}$, $a \in [c]$, must be different. In the simplification from (A.8) to (A.11) we also replaced the permutation $\nu_\sigma$ by a new permutation $\mu$. For this step, recall that $\nu_\sigma$ in (A.8) is the permutation that brings $h_{\sigma(c+1)}, \ldots, h_{\sigma(m)}$ in increasing order with $h_{\sigma(c+1)} = k_{c+1} < h_{\sigma(c+2)} = k_{c+2} < \cdots < h_{\sigma(c+p)} = k_{c+p}$, which implies that $\nu_\sigma(j) = j$ for all $j \in \{c+1, \ldots, c+p\}$. Thus the sign of $\nu_\sigma$ is equal to the sign of $\nu_\sigma|_{\{c+p+1,\ldots,m\}}$. The latter restriction is denoted by $\mu$ in (A.11).

Noting that $k_b \in \bar{J} \cap \bar{K}$ if $b \in \{c+p+1, \ldots, m\}$, we see that

$$\sum_{\mu \in \mathcal{S}_{\bar{I} \cap \bar{L}}} (-1)^\mu \prod_{b=c+p+1}^{m} t_{k_b \mu(b)} = \det(T_{(\bar{J} \cap \bar{K}) \times (\bar{I} \cap \bar{L})}).$$

Thus, we have shown that $\mathrm{E}[\det(W_{I \times J}) \det(W_{K \times L})]$ is equal to the expectation of

$$(A.12) \quad \left(\prod_{i \in I} t_{ii}^2\right) \det(T_{(\bar{L} \cap \bar{I}) \times (\bar{J} \cap \bar{K})})^2 \det(T_{(\bar{L} \cap \bar{J}) \times (\bar{I} \cap \bar{K})})^2 \left[\prod_{a=1}^{c} \left(\sum_{h=a}^{[k_a] \setminus \bar{I}} t_{k_a h}^2\right)\right].$$

Since $t_{ii}^2 \sim \chi_{n-i+1}^2$, and moreover,

$$\sum_{h=a}^{[k_a] \setminus \bar{I}} t_{k_a h}^2 \sim \chi_{(n-k_a+1)+(k_a-a)-|\bar{I}|}^2 = \chi_{(n-a+1)-(m-c)}^2 = \chi_{n-m+c-a+1}^2,$$

this proof can be completed using the results on expected values from the proof of Theorem 4.4.

**A.3. A noncentral Wishart determinant.** As in Lemma 5.4, we consider $X \in \mathbb{R}^{m \times m}$ distributed according to $\mathcal{N}_{m^2}(A, I_m \otimes I_m)$, $A = (a_{ij}) \in \mathbb{R}^{m \times m}$. From the independence of the entries of $X$, it follows that

$$(A.13) \quad \mathrm{E}[\det(X)] = \det(A).$$

If $A$ is nonzero, then the matrix $XX^T$ follows a *noncentral Wishart distribution*. Theorem 10.3.7 in [12] provides a general formula for moments of the determinant of a noncentral Wishart matrix in terms of hypergeometric functions with matrix argument. Here, $X$ is a square matrix and we can give a simple formula $\mathrm{E}[\det(X)^2] = \mathrm{E}[\det(XX^T)]$ that involves only traces and compounds.

LEMMA A.1. *The expectation of* $\det(XX^T) = \det(X)^2$ *can be expressed as*

$$\mathrm{E}[\det(X)^2] = \sum_{k=0}^{m} (m-k)! \cdot \mathrm{tr}[(AA^T)^{(k)}].$$

Here, $(AA^T)^{(0)} := 1 \in \mathbb{R}$ and $(AA^T)^{(m)} = \det(AA^T)$. By (A.13),

$$\text{Var}[\det(X)] = \sum_{k=0}^{m-1}(m-k)! \cdot \text{tr}[(AA^T)^{(k)}].$$

PROOF. Let $\mathcal{S}_m$ be the group of permutations of $[m]$. Then,

$$\text{E}[\det(X)^2] = \sum_{\sigma \in \mathcal{S}_m} \sum_{\tau \in \mathcal{S}_m} \prod_{j=1}^{m} \text{E}[X_{\sigma(j)j} X_{\tau(j)j}]$$

$$= \sum_{\sigma \in \mathcal{S}_m} \sum_{\tau \in \mathcal{S}_m} (-1)^{\sigma+\tau} \prod_{j=1}^{m}(\delta_{\sigma(j)\tau(j)} + a_{\sigma(j)j}a_{\tau(j)j}),$$

where $\delta_{ij}$ is the Kronecker delta. The product

$$\prod_{j=1}^{m}(\delta_{\sigma(j)\tau(j)} + a_{\sigma(j)j}a_{\tau(j)j}) = \sum_{J \subseteq [m]} \left(\prod_{j \in J} a_{\sigma(j)j}a_{\tau(j)j}\right) \cdot \left(\prod_{j \notin J} \delta_{\sigma(j)\tau(j)}\right).$$

Therefore, if we define $g_J(\sigma) = \sum_{\substack{\tau \in \mathcal{S}_m \\ \tau(j)=\sigma(j) \forall j \notin J}} (-1)^{\sigma+\tau} \prod_{j \in J} a_{\sigma(j)j}a_{\tau(j)j}$, then

$$\text{E}[\det(X)^2] = \sum_{k=0}^{m} \sum_{J \in \{^m_k\}} \sum_{\sigma \in \mathcal{S}_m} g_J(\sigma).$$

Note that the permutations appearing in the definition of $g_J(\sigma)$ satisfy $\tau(J) = \sigma(J)$.

Let $\sigma_1, \sigma_2 \in \mathcal{S}_m$ be two permutations such that $\sigma_1(j) = \sigma_2(j)$ for all $j \in J$. Moreover, let $\tau_1, \tau_2 \in \mathcal{S}_m$ satisfy $\tau_1(j) = \tau_2(j)$ for all $j \in J$, $\tau_1(j) = \sigma_1(j)$ for all $j \notin J$, and $\tau_2(j) = \sigma_2(j)$ for all $j \notin J$. Then it holds for the permutation signs that $(-1)^{\sigma_1}(-1)^{\tau_1} = (-1)^{\sigma_2}(-1)^{\tau_2}$. This implies that $g_J(\sigma_1) = g_J(\sigma_2)$. We obtain

$$\text{E}[\det(X)^2]$$
(A.14)
$$= \sum_{k=0}^{m}(m-k)! \sum_{J \in \{^m_k\}} \sum_{I \in \{^m_k\}} \sum_{\bar{\sigma} \in \mathcal{S}_k} \sum_{\bar{\tau} \in \mathcal{S}_k} (-1)^{\bar{\sigma}+\bar{\tau}} \prod_{h=1}^{k} a_{i_{\bar{\sigma}(h)}j_h} a_{i_{\bar{\tau}(h)}j_h}.$$

By the Binet–Cauchy theorem,

$$\text{E}[\det(X)^2] = \sum_{k=0}^{m}(m-k)! \sum_{I \in \{^m_k\}} \sum_{J \in \{^m_k\}} \det(A_{IJ})^2$$

$$= \sum_{k=0}^{m}(m-k)! \sum_{I \in \{^m_k\}} \det(A_{I \times [m]} A^T_{I \times [m]})$$

$$= \sum_{k=0}^{m}(m-k)! \cdot \text{tr}[(AA^T)^{(k)}]. \qquad \square$$

**Acknowledgments.** We would like to thank the Editor and two referees for helpful comments.

M. DRTON
DEPARTMENT OF STATISTICS
UNIVERSITY OF CHICAGO
CHICAGO, ILLINOIS 60637
USA
E-MAIL: drton@galton.uchicago.edu

H. MASSAM
DEPARTMENT OF MATHEMATICS AND STATISTICS
YORK UNIVERSITY
TORONTO, ONTARIO
CANADA M3J 1P3
E-MAIL: massamh@yorku.ca

I. Olkin
Department of Statistics
Stanford University
Stanford, California 94305-4065
USA
E-mail: iolkin@stat.stanford.edu